\documentclass{article}
\usepackage{amsmath}
\usepackage{amsfonts} 
\usepackage{amssymb}
\usepackage{hyperref}
\usepackage{diagrams}
\diagramstyle[PostScript=dvips]
\usepackage{amsthm}
\newarrow{Dotsto}....>
\newarrow{Equals}=====
\def\leftsimeq{\!\wr\!\!\mid}
\def\rightsimeq{\mid\wr}
\def\simeqq{\overset{\sim}{\relbar}}
\def\backsimeqq{\overset{\backsim}{\relbar}}
\newarrowmiddle{simeq}\simeqq\backsimeqq\leftsimeq\rightsimeq
\newarrow{Isom}--{simeq}->
\newarrow{Isomto}--{simeq}->

\newcommand{\myqed}{\nobreak \ifvmode \relax \else
      \ifdim\lastskip<1.5em \hskip-\lastskip
      \hskip1.5em plus0em minus0.5em \fi \nobreak \hfill
      \vrule height0.5em width0.5em depth0.25em\fi}

\newcommand{\1}{1\!\!1}
\newcommand{\Grobner}{ Gr\"{o}bner }
\newcommand\Poincare{{Poincar\'{e} }}
\newtheorem{Theorem}{Theorem}
\numberwithin{Theorem}{section}
\newtheorem{Lemma}[Theorem]{Lemma}
\newtheorem{Proposition}[Theorem]{Proposition}
\theoremstyle{definition}
\newtheorem{Definition}[Theorem]{Definition}
\newtheorem{LemmaDefinition}[Theorem]{Lemma/Definition}
\newtheorem{Example}[Theorem]{Example}
\newtheorem{Remark}[Theorem]{Remark}
\newcommand{\modulo}[1]{{ \ \ \left(\mathbf{mod} \ {#1}\right)}}
\renewcommand{\deg}{\mathbf{deg}}
\DeclareMathOperator{\tensor}{\otimes}
\DeclareMathOperator{\comp}{\circ}
\DeclareMathOperator{\len}{\mathbf{len}}
\newcommand{\g}{\mathfrak{g}}
\newcommand{\h}{\mathfrak{h}}
\newcommand{\GG}{\mathbb{G}}
\newcommand{\HH}{\mathbb{H}}
\newcommand{\ZZ}{\mathbb{Z}}
\newcommand{\QQ}{\mathbb{Q}}
\newcommand{\FF}{\mathbb{F}}
\newcommand{\T}{\mathbf{T}}

\newcommand{\U}{\mathfrak{U}}
\newcommand{\uvec}[1]{{\underline{{#1}}}}
\newcommand{\Dist}{\mathbf{Dist}}
\newcommand{\Comm}{\mathbf{Comm}}
\newcommand{\Hom}{\mathbf{Hom}}

\newcommand{\requireEquiv}{{\stackrel{\mbox{\tiny{REQ}}}{\equiv}}}
\DeclareMathOperator{\inv}{\mathbf{inv}}
\DeclareMathOperator{\LM}{\mathbf{LM}}
\DeclareMathOperator{\LC}{\mathbf{LC}}
\DeclareMathOperator{\gge}{\uvec{\gg}}

\newcommand{\canon}{{\mathbf{c}}}

\newcommand\Ref[2]{(#1. \ref{#2})}

\def\enumEnvironment#1#2{
  \newenvironment{#1}
  {
    \begin{#2}
      {\hskip1pt }\\{ \hskip1pt}
      \vspace{-1.6em}
      \numberwithin{enumi}{Theorem}
      \begin{enumerate}
        \setlength{\labelwidth}{0em}
        \setlength{\labelsep}{0.5em}
        \setlength{\leftmargin}{0em}
        \setlength{\itemindent}{0.6em}
        \numberwithin{enumi}{Theorem}
        \renewcommand{\theenumi}{{\bf{\theTheorem.\arabic{enumi}}}}
      }
      {
      \end{enumerate}
    \end{#2}
  }
}

\def\makeEnumeratedEnvironment#1
{
  \expandafter\message{Creating environment e#1.}
  \expandafter\enumEnvironment{e#1}{#1}
}

\makeEnumeratedEnvironment{Definition}
\makeEnumeratedEnvironment{Remark}
\makeEnumeratedEnvironment{Example}

\author{Frederick Leitner}
\title{Non-Commutative Formal Groups in Positive Characteristic}

\begin{document}

\maketitle

\begin{abstract}
  \noindent We describe geometric non-commutative formal groups
  in terms of a geometric commutative formal group with a Poisson
  structure on its splay algebra.  We describe certain natural properties
  of such Poisson structures and show that any such Poisson structure
  gives rise to a non-commutative formal group.
\end{abstract}

\noindent We describe geometric non-commutative formal groups
in terms of a geometric commutative formal group with a Poisson
structure on its splay algebra.  We describe certain natural properties
of such Poisson structures and show that any such Poisson structure
gives rise to a non-commutative formal group.

In  characterizing non-commutative formal groups $\GG$, it is 
natural to consider their distribution algebras, $\Dist(\GG)$. Although we
are interested mainly in the characteristic $p>0$ setting, 
we make an effort to work by analogy to the characteristic zero version of
our problem.   Let $\GG$ be a formal group 
with Lie algebra $\g$.  We recall that, in characteristic zero, 
we have a canonical isomorphism between $\Dist(\GG)$ and the universal 
enveloping algebra $\U\g$.  Then our problem reduces to describing
the algebraic relationship between the symmetric $\mathbf{S}\g$, $\U\g$, and
the tensor algebra $\T\g$ as in the following diagram:
\begin{equation}\label{triangle}
  \begin{diagram}[height=1.7em]
    & &\T\g && \\
    &\ldTo & &\rdTo&   \\
    \mathbf{S}\g& &\rIsomto_{\mathbf{PBW}}& & \U\g.
  \end{diagram}
\end{equation}
Here, the diagonal arrows are algebra maps and the horizontal arrow is  a
vector space  isomorphism given by the 
\Poincare\!\!-Birkhoff-Witt (PBW) theorem.

In the positive characteristic setting,  there are a few features that
present difficulties which are absent in the zero characteristic setting.
In a sense made precise
below, the problem is non-linear and, unlike in the zero characteristic
setting, the algebras of importance are not quadratic.  

Let us recall that in characteristic zero, the PBW theorem gives an identification between 
$\mathbf{S}\g$ and the associated graded algebra of $\U\g$.
We may  regard $\mathbf{S}\g$ as the universal
enveloping algebra for the trivial Lie algebra structure on $\g$, so that $\mathbf{S}\g$ may
be viewed as the distribution algebra of the unique commutative formal group of dimension
$n=\dim(\g)$.  Indeed, $\mathbf{S}\g=\Dist(\GG_a^n)$ where $\GG_a$ is the additive formal group.
In fact, one may identify, regardless of the characteristic, the associated graded algebra
of $\Dist(\GG)$ with $\Dist(\GG_a^n)$.
However, in positive characterisitc  the difficulty  arises
there are many choices of commutative formal groups of a fixed dimension
(\cite{Dieudonne},  \cite{Lazard}, \cite{Manin}).  Thus,
one might suspect that passing to the associated graded algebra is too
course as it loses too much of the underlying commutative geometry of the formal group.
It is for these reasons that we introduce geometric formal groups and 
describe Poisson structure which preserve these commutative sub-structures.  

In  $\S1$, we introduce the notion of ``splay'' algebras and prove
a PBW type theorem for them.  In $\S2$, we recall some relevant
details of distribution algebras of formal groups.  In  $\S3$, we characterize completely
geometric non-commutative formal groups in terms of certain Poisson structures
on their splays.

We enforce several assumptions and notations throughout this paper.
By $R$ we will denote a commutative ring with unit, and if it is a field,
we will denote it by $k=R$.
All algebras are associate unital $R$-algebras, unless
they are Lie algebras. All tensor products will be
taken over $R$.  If $V$ is a free $R$-module,
we denote by $\T V$ the tensor algebra on $V$ (over $R$), and $\mathbf{S} V$ the symmetric
algebra on $V$.  If $\mathcal{X}$ is a set, we denote by
$\T\mathcal{X}$  the free  algebra over $R$ generated by $\mathcal{X}$ and
by $\mathbf{S}\mathcal{X}$ the free commutative algebra generated by $\mathcal{X}$.
All ideals will be two-sided ideals.  
If $A$ and $B$ are algebras, we denote by $A\tensor B\rTo^{\tau} B\tensor A$
the {\bf transposition} or {\bf swapping} map  $a\tensor b \rMapsto^{\tau}
b\tensor a$.
We assume that none of the algebras $A$ under consideration have
quasi-zeros, i.e. elements $a\in A$ such that for all $b,c\in A$ one has $bac=0$.
We also  make implicit use of the theory
of non-commutative \Grobner basis, for which we follow the notation and
conventions of \cite{Li} in the case of a field $k$ and for general $R$
those in \cite{LeiPaw}.
By abuse of language, non-commutative often
means not necessarily commutative.

I am indebted to  P. Bressler, P. Foth, M. Kim, K. Lux, D. Ulmer, and A. Vasiu
for input and conversations.  During the  time of preparation, I was
partially supported by a  NSF VIGRE grant while at the University of Arizona,
and by a VATAT fellowship while at the Ben-Gurion University of the Negev.

\section{Poisson Algebras with Internal Symmetry}\label{PoissonAlgebras}

In this section, we generalize  the algebraic relationship
as indicated in \Ref{Eqn}{triangle} for a Lie algebra $\g$ over $R$.
Choosing any $R$-module basis  $\g=R\langle X_i
\rangle_{i\in I}$, one can view $\mathbf{S}\g$ as the tensor product of the one dimensional
algebras $\mathbf{S}\{X_i\}$.  It is from this vantage point that we wish to
generalize.  We will be looking at $S$, a commutative algebra, which has a
fixed identification as the tensor product $S=\bigotimes_{i\in I} A_i$ of
commutative sub-algebras $A_i$ for $i$ in some indexing set $I$.  We then
desire a non-commutative algebra $U$ for which there is a  PBW-type isomorphism between $S$
and $U$, and an algebra $T$ for which both $U$ and $S$ are quotients by
appropriate commutator relations:
\begin{displaymath}
  \begin{diagram}[height=1.7em]
    &&T && \\
    & \ldTo& &\rdTo& \\
    S & &\rIsomto_{\mathbf{PBW}} &  & U.
  \end{diagram}
\end{displaymath}
Specifically, we desire that $S$ be the quotient of $T$ by ideal of
relations generated by $fg-gf$ for $f,g\in T$.
We recall that  a Poisson structure on an algebra $A$ is a
skew-symmetric bi-derivation $\pi$  which  satisfies the Jacobi identity.

\begin{eDefinition}
\item
  Let $I$ be an indexing set and for each $i$ fix a commutative algebra
  $A_i$.
  We define the   {\bf splay} of the $A_i$ to be
  the following diagram $\overline{T}$ in the category of algebras:
  \begin{displaymath}
    A_i \rInto^{\iota_i} \coprod_{i\in I} A_i.
  \end{displaymath}
  We call $T:=\coprod_i A_i$, the co-product of the $A_i$,
  the {\bf splay algebra} of $\overline{T}$
  and we call the $A_i$ the {\bf coordinate algebras} of $\overline{T}$.
\item
  Let $\overline{T}$ (resp. $\overline{T}'$) be a  splay
  with coordinate algebras $A_i$ for $i\in I$ (resp. $B_j$ for $j\in J$).
  A {\bf morphism} of splays $\overline{T} \rTo^{\Phi} \overline{T}'$
  is an algebra morphism $T\rTo^{\Phi}T'$ of the respective splay algebras
  such that if for every $i \in I$ there is some $j_i \in J$ and
  some $A_i \rTo^{\Phi_{j_i}}B_j$ such that
  $\Phi$ is the map induced from the UMP of the co-product $T$ as follows:
  \begin{displaymath}
    \begin{diagram}[height=1.7em]
      A_i & \rInto^{\iota_i} & \coprod A_i \\
      \dTo^{\Phi_{j_i}} &\ldDashto(2,4)^{\exists !} &  \\
      B_{j_i} & & \\
      \dTo^{\iota_{j_i}} & & \\
      \coprod_i B_j
    \end{diagram}
  \end{displaymath}
\item
  Let $\pi$ be a Poisson structure on an algebra $A$, and suppose that
  the $A_i$ are sub-algebras.
  Then we say that $\pi$ has {\bf internal symmetry}
  (with respect to the $A_i$) if  $\pi(A_i,A_i) =0$.
  If  $\overline{T}$ is a splay with coordinate algebras $A_i$, then  a Poisson
  structure  $\pi$ on $\overline{T}$ is a Poisson structure $\pi$ on
  $T$ with internal symmetry with respect to the $A_i$'s.
\item
  Let $T$ and $T'$ be splay algebras equipped with 
  Poisson structures $\pi$ and $\pi'$ respectively.  We say that
  a morphism  of splays $\overline{T}\rTo^{\Phi} \overline{T}'$  is a {\bf Poisson morphism}
  if $\Phi\comp \pi = \pi' \comp (\Phi\tensor \Phi)$.  If $T=T'$ as splays, 
  then we say  that $\pi$ and $\pi'$ are {\bf equivalent} Poisson
  structures and that  $\Phi$ gives an {\bf equivalence} between $\pi$
  and $\pi'$.
\item
  Let $A$ be an algebra and $\pi$ a Poisson structure on $A$.  We say that
  $\pi$ {\bf vanishes on constants} if for all $a\in A$ we have
  $\pi(1,a)=0=\pi(a,1)$.
\item
  Suppose that $A=\bigcup_{n\ge 0} A^{(n)}$ is a filtered algebra
  and $\pi$ is a Poisson structure on $A$.  We say that
  $\pi$ is {\bf strongly filtered} if:
  \begin{displaymath}
    \pi(A^{(m)},A^{(n)}) \subseteq A^{(m+n-1)}.
  \end{displaymath}
\end{eDefinition}

\noindent We will assume for the remainder of this paper
that all Poisson structures vanish on constants.
\begin{Example}
  For an algebra $A$  one has the
  {\bf canonical Poisson   structure} defined by:
  \begin{displaymath}
    \pi_{\canon}(a,b):= ab-ba
  \end{displaymath}
  for $a,b\in A$.  The fact that it is a Poisson structure follows
  from the associativity of the algebra $A$.  As we have assumed that
  $A$ has no quasi-zeros, clearly $\pi_{\canon}$
  vanishes if and only if $A$ is a commutative algebra, so that
  $\pi_{\canon}$ is a measure of the non-commutativity of $A$.
  \myqed
\end{Example}

\begin{Remark}
  Looking ahead, the Poisson structures $\pi$ on the splay algebras  measure
  the non-commutativity of the algebras $U$.   In our motivating example of $\g$ 
  one  takes $T=\T\g$, $S=\mathbf{S}\g$, and $U=\U\g$.  One usually considers the 
  Kostant-Kirilov-Souriau $\pi_{\scriptscriptstyle{\mathbf{KKS}}}$  Poisson
  bracket which is a Poisson bracket on $\mathbf{S}\g$ rather than $\T\g$ (c.f. \cite{KorSoi}).  
  However, one may may also form a Poisson bracket $\tilde{\pi}$ on $\T\g$ by extending the Lie bracket by the
  bi-derivation law.   Working modulo the commutator relations  $\eta\zeta - \zeta\eta$ for $\eta,\zeta\in \g$ we have
  for all $f\in\T\g$ that:
  \begin{eqnarray*}
    \tilde{\pi}(\eta\zeta -\zeta\eta,f) &\equiv& \eta\tilde{\pi}(\zeta,f) + \tilde{\pi}(\eta,f)\zeta \\
    & &  - \zeta\tilde{\pi}(\eta,f) - \tilde{\pi}(\zeta,f)\eta \\
    &\equiv & 0
  \end{eqnarray*}
  so that $\tilde{\pi}$ on $\T\g$ descends to $\pi_{\scriptscriptstyle{\mathbf{KKS}}}$ on $\mathbf{S}\g$.

  We cannot hope to have such a nice situation in general.  For example, suppose that
  $R$ is a $\ZZ/m\ZZ$-algebra for some $m\ge 2$, and suppose that we have a relation of the
  form $f^m=g$ in $S$. Then in order for a Poisson structure $\pi$ on $T$ to descend to $S$
  we would need for all $h\in T$ that:
  \begin{eqnarray*}
    0&\requireEquiv & \pi(f^m-g,h) \\
    &\equiv &m! f^{m-1}\pi(f,h) - \pi(g,h) \\
    &\equiv & -\pi(g,h)
  \end{eqnarray*}  
  so that we must have  either $g=0$ or $\pi(g,h)=0$.  Although this will not happen in general,
  we will see below that  $S=\Dist(\GG_a^n)$,  over any $\FF_p$-algebra  $R$,
  falls into the former case, and thus a Poisson structure on the splay of
  some $\Dist(\GG_a)$'s
  will descend.  
  \myqed
\end{Remark}

\noindent There are two main types of PBW proofs.  One makes use
of the theory of \Grobner bases as in \cite{Mora} and \cite{deGraaf}.  The
other is from the viewpoint of deformations  of Koszul algebras as in
\cite{BravGaits} and \cite{floy}.  Below, we generalize the
\Grobner basis setting.

\begin{eRemark}
\item
  Let $\mathcal{X}$ be a set and suppose that there is a map
  $ \mathcal{X} \rTo^{\deg} \ZZ_{\ge 1}$. Suppose that $<$ is a total
  ordering on the set $\mathcal{X}$.  Then if $A$ is either
  $\mathbf{S}\mathcal{X}$ or $\T\mathcal{X}$, we obtain a monomial ordering
  in the sense of \cite{Li} and \cite{LeiPaw}  on $A$ by using the graded lexicographic
  ordering where for a monomial $\uvec{m}=x_1\cdots x_n$ with $x_i\in \mathcal{X}$,
  we define its degree and length by: 
  \begin{eqnarray*}
    \deg(\uvec{m})&:=&\sum_{1\le i  \le n} \deg(x_i) \\
    \len(\uvec{m})&:=&n.
  \end{eqnarray*}
  We define:
  \begin{displaymath}
    A^{(n)}:=R\langle \uvec{m}  \ \ |  \ \deg(\uvec{m})\le n \rangle \ \  \subseteq A
  \end{displaymath}
  to be the $R$-span of monomials of degree at most $n$, we see that $A$ is a
  filtered algebra.  $A$ is  in fact graded algebra where:
  \begin{eqnarray*}
    A^{n}&:=&R\langle \uvec{m}  \ \ |  \ \deg(\uvec{m})=n \rangle \ \
    \subseteq A \\
    A &=& \bigoplus_{n \ge 0} A^{n}.
  \end{eqnarray*}
  If, in particular, we
  take $\mathcal{X} \rTo^{\deg} \ZZ_{\ge 1}$ to be the constant function
  $1$, then $\deg(\uvec{m})=\len(\uvec{m})$ and
  we recover the usual grading
  on $\mathbf{S}\mathcal{X}$ and   $\T\mathcal{X}$.  
  If $\mathcal{J}$ is an ideal in $A$ and $B:=A/\mathcal{J}$, then $B$
  inherits  a filtration from $A$:
  \begin{displaymath}
    B=\bigcup_{n\ge 0} B^{(n)} \ \ \ \ \ \ \ B^{(n)}:=q(A^{(n)})
  \end{displaymath}
  where $A\rTo^{q} B$ is the quotient map.
  We will assume below that all algebras constructed in this manner are
  equipped with this filtration.
\item\label{embedding}
  Let $I$ be an indexing set, and for each $i\in I$ suppose we have sets
  $\mathcal{X}_i$, and suppose that for each $i \in I$ we have an ordering
  $<_i$ on $\mathcal{X}_i$.  
  Define $\mathcal{X}:=\coprod_{i\in I} \mathcal{X}_i$.
  Then we have an embedding:
  \begin{equation}\label{embeddingEqn}
    \mathbf{S}\mathcal{X}_i \rInto^{\iota_i} \T\mathcal{X}_i \rInto  \T\mathcal{X}
  \end{equation}
  given by writing any  monomial $\uvec{m}$ of length $n$ in $\mathcal{X}_i$ as:
  \begin{displaymath}
    \uvec{m}=x_1 x_2 x_3 \cdots x_n
  \end{displaymath}
  where $x_i \in \mathcal{X}_i$ and $x_1 \le_i x_2 \le_i \cdots \le_i x_n$,
  and we define:
  \begin{displaymath}
    \iota_{i}(\uvec{m}) := x_1 x_2 x_3 \cdots x_n \in \T\mathcal{X}_i.
  \end{displaymath}
  The map from $\mathbf{S}\mathcal{X}_i$ is then the  $R$-linear extension of this map.
  By abuse of notation, we denote $\iota_i(\uvec{m})$ simply by $\uvec{m}$
  and the composite map of \Ref{Eqn}{embeddingEqn} also by $\iota_i$.
\end{eRemark}

\begin{Lemma}\label{Tbasis}
  Let $\overline{T}$ be a splay with coordinate algebras $A_i$ over a field $k$.
  Suppose that each $A_i$ has a fixed  presentations:
  \begin{displaymath}
    \mathcal{S}_i \rInto \mathbf{S}\mathcal{X}_i \rOnto A_i
  \end{displaymath}
  for some set $\mathcal{X}_i$ and ideal of relations $\mathcal{S}_i$.
  Suppose that for each $i$ there is
  a map:
  \begin{displaymath}
    \deg_i : \mathcal{X}_i \rTo \ZZ_{\ge 1}
  \end{displaymath}
  and an ordering $<_i$ on $\mathcal{X}_i$ such that
  for some indexing set   $J_i$ there
  is a collection of polynomials $\{f_{i,j} \ | \ j \in J_i \}$
  which forms a \Grobner basis for the ideal $\mathcal{S}_i$ with
  respect to the graded lexicographic ordering on $\mathbf{S}\mathcal{X}_i$ induced
  by $<_i$ and $\deg_i$.
  Then the splay algebra $T$ has a $k$ basis of elements of the
  form:
  \begin{displaymath}
    \uvec{a}_1 \uvec{a}_2 \cdots \uvec{a}_k
  \end{displaymath}
  where for each $1\le j \le k$ there is some $i_j \in I$  
  such that   $i_j \ne i_{j+1}$ and such that
  $\uvec{a}_j=a_{j,1}\cdots a_{j,n_j}$
  are non-decreasing monomials of length $n_j\ge 1$ in
  $\mathcal{X}_{i_j}$ which 
  are not leading monomial ideal of the ideal $\mathcal{S}_{i_j}$, 

\end{Lemma}

\begin{proof}
  Adopting the notation of \Ref{Rmk}{embedding}, 
  one may exhibit $\coprod_i A_i$ as:
  \begin{displaymath}
    T = \T\mathcal{X}/\mathcal{R}
  \end{displaymath}
  where $\mathcal{R} = \sum_i \mathcal{R}_i$ and $\mathcal{R}_i$ is the
  ideal generated by $\mathcal{S}_i$ and the relations $s_{a,a'}:=aa'-a'a$ for
  $a,a'\in A_i$.

  The statement now reduces to showing that the set:
  \begin{displaymath}
    \bigcup_{i\in I} \{ f_{i,j} \ | \ j\in J_i\} \cup \bigcup_{i\in I} \{ s_{a,a'} \ |
    \     a>_i a' \in \mathcal{X}_i \}
  \end{displaymath}
  is a \Grobner basis for the ideal $\mathcal{R}$ for some monomial ordering
  on $\T\mathcal{X}$.
  This follows readily
  from the  graded lexicographic ordering given by $\deg$ and
  an ordering $\ll$ as follows.
  As $\mathcal{X}$ is the co-product in sets of the $\mathcal{X}_i$,
  the maps $\deg_{i}$ induce a map $\mathcal{X}\rTo^{\deg} \ZZ_{\ge 1}$.
  In concrete terms, if $\uvec{m}=\uvec{m_1}\cdots \uvec{m_n}$ is a monomial in $\T\mathcal{X}$
  with the $\uvec{m_j}$ monomials in $\mathcal{X}_{i_j}$ for some $i_j \in I$, then we have:
  \begin{displaymath}
    \deg(\uvec{m}):= \sum_{1\le j \le n}  \deg_{i_j}(\uvec{m_j}).
  \end{displaymath}
  Let $x_i\in \mathcal{X}_i \subseteq \mathcal{X}$  and $x_{i'} \in
  \mathcal{X}_{i'}\subseteq \mathcal{X}$ we define:
  \begin{displaymath}
    x_i \ll x_{i'} \ \ \mbox{ if } \ \
    \left\{
      \begin{array}{l}
        i < i', \\
        \mbox{OR} \\
        i=i' \ \mbox{ and } \  x_i <_i x_{i'}.
      \end{array}
    \right.
  \end{displaymath}
\end{proof}

\begin{eRemark}
\item
  By the results in \cite{LeiPaw},
  we may easily extend \Ref{Lem}{Tbasis} to the case of a general ring $R$
  as long as we assume that the leading coefficients of the $f_{i,j}$
  are units of $R$.
\item
  \label{embedding2}
  Under the assumptions of \Ref{Lem}{Tbasis} and the ordering $\ll$
  introduced in its proof,  we obtain an  $R$-linear map:
  \begin{displaymath}
   \begin{diagram}
      T &  \rInto^{\iota_{\ll}}  & \T\mathcal{X}
    \end{diagram}
  \end{displaymath}
  such that $q\comp \iota_{\ll}$ is the identity map on $T$ where $q$
  is the quotient map $\T\mathcal{X}\rOnto^{q} T$.
  If $\mathcal{J}$ is an ideal of $T$, we denote by $\tilde{J}_{\ll}$
  the ideal of $T\mathcal{X}$ generated by image of
  $\iota_{\ll}(\mathcal{J})$.

  Suppose now that $\pi$ is a Poisson structure on a splay $\overline{T}$.
  Then we obtain a Poisson structure on $\T\mathcal{X}$, denoted by
  $\tilde{\pi}$, 
  by applying the bi-derivation law to the
  following map:
  \begin{displaymath}
    \begin{diagram}[height=1.7em]
      \mathcal{X}\times\mathcal{X}&\rTo& A&\rInto^{\iota_{\ll}}&\T\mathcal{X}\\
      (x_i,x_{i'}) &\rMapsto &  \pi(x_i,x_i') &\rMapsto & \iota_{\ll}\comp \pi(x_i,x_{i'}).
    \end{diagram}
  \end{displaymath}
\end{eRemark}

\begin{Lemma}\label{measure}
  Let $\overline{T}$ be a splay over a field $k$ with coordinate algebras $A_i$ 
  which  satisfy the conditions of the \Ref{Lem}{Tbasis}.
  Suppose that $\pi$  Poisson structure on $\overline{T}$
  Denote by $\mathcal{J}$ the  $T$  ideal:
  \begin{displaymath}
    \mathcal{J} = \left( a a' - a' a - \pi(a,a') \right).
  \end{displaymath}
  Let $\tilde{\pi}_{\ll}$ and  $\tilde{\mathcal{J}}_{\ll}$
  be as in   \Ref{Rmk}{embedding2} and $\mathcal{R}$ be as in the
  proof of \Ref{Lem}{Tbasis}.
  Then:
  \begin{itemize}
  \item 
    The Poisson structure 
    $\tilde{\pi}_{\ll}$ descends to $\pi$.
  \item
    Suppose that $f,g\in\T\mathcal{X}$, then 
    \begin{displaymath}
      fg-gf \equiv \tilde{\pi}_{\ll}(f,g) \modulo{\mathcal{R} +
      \tilde{\mathcal{J}}_{\ll}}.
    \end{displaymath} 
  \item 
    $\pi$ descends to the canonical Poisson structure
    on $U:=T/\mathcal{J}$.
  \end{itemize}
\end{Lemma}
\begin{proof}
  The first statement is clear by the construction of $\tilde{\pi}_{\ll}$.
  The third statement follows readily from the second. For the second, by
  $k$-linearity,
  we may assume that $f$ and $g$ are monomials.  In which
  case, we proceed by induction on the sum of their lengths $L:=\len(f) +
  \len(g)$. For $L \le 2$ it is trivial by definition of $\tilde{\pi}_{\ll}$.
  Suppose now that we have proved it for all $L' \le L$.  Let $x\in
  \mathcal{X}$
  and consider $f'=fx$ where $\len(f)\le L'$.  Then, we have:
  \begin{eqnarray*}
    f'g-gf' &\equiv& fxg-gfx \\
    &\equiv& (f\tilde{\pi}_{\ll}(x,g) + fgx) -(fgx+ \tilde{\pi}_{\ll}(g,f)x) \\
    &\equiv& f\tilde{\pi}_{\ll}(x,g) + \tilde{\pi}_{\ll}(f,g)x  \\
    &\equiv& \tilde{\pi}_{\ll}(fx,g) \equiv \tilde{\pi}_{\ll}(f',g).
  \end{eqnarray*}
  A similar statement for $\tilde{\pi}_{\ll}(g,f')$ gives the desired result.
\end{proof}

\begin{Remark}
  Let $\overline{T}$ (resp. $\overline{T}'$) be a splay with Poisson
  structure $\pi$ (resp. $\pi'$) and suppose that $\overline{T}\rTo^{\Phi}
  \overline{T}'$ is a Poisson morphism.  Let $U$ (resp. $U$) be as in
  \Ref{Lem}{measure}.  Then it  is clear that $\Phi$ descends
  to an algebra morphism $U\rTo^{\Phi} U'$.  If $A_i$ for $i\in I$
  (resp. $A'_j$ for $j\in J$) are the
  coordinate algebras of $\overline{T}$ (resp. $\overline{T}'$), then we
  have the following commutative diagram:
  \begin{displaymath}
    \begin{diagram}[height=1.7em]
      A_i & \rInto^{\iota_i}  & T & \rOnto &U \\
      \dTo^{\Phi_{j_i}} & & & & \dTo_{\Phi} \\
      A'_{j_i} & \rInto^{\iota_{j_i}}  & T' & \rOnto &U'.
    \end{diagram}
  \end{displaymath}
  \myqed
\end{Remark}

\noindent 
We now prove a PBW-type theorem.  In the case of a Lie algebra $\g$ this reduces
to the usual PBW theorem where one has chosen an ordered basis for
the Lie algebra (c.f. \cite{Serre}, \cite{Mora}, \cite{deGraaf}).
Again, we may extend \Ref{Lem}{measure} and \Ref{Rmk}{PBW} to the case of a
general ground ring $R$ as long as the leading coefficients of the $f_{i,j}$'s are units.

\begin{Theorem}\label{PBW}
  Suppose $\overline{T}$ is a splay over a field $k$ with coordinate algebras $A_i$
  which  satisfy the conditions of \Ref{Lem}{Tbasis}.
  Suppose that $\pi$
  is a strongly filtered Poisson structure on $\overline{T}$,
  and let $\mathcal{J}$ be as in \Ref{Lem}{measure}.
  Then, fixing the choice of a total ordering $<$ on $I$, 
  $U$ is isomorphic as a $k$-module to $S:=\bigotimes_i A_i$.
\end{Theorem}

\begin{proof}
  It is now enough to show that the set:
  \begin{displaymath}
   \mathcal{Y}:= \{ f_{i,j} \ | i\in I, j\in J_i \} \cup \{ g_{a,b} \ | \ a \gg b \in
    \mathcal{X} \} 
  \end{displaymath}
  is a \Grobner basis  of $\mathcal{R} + \tilde{\mathcal{J}}_{\ll}$ with
  respect to $\deg$ and the ordering $\ll$ as in the proof of \Ref{Lem}{Tbasis}.

  We first show that $\mathcal{Y}$ generates $\mathcal{R}+\tilde{\mathcal{J}}_{\ll}$.
  By
  $k$-linearity of the Poisson bracket, and as $\mathcal{R}=\sum_i
  \mathcal{R}_i$  where $\mathcal{R}_i\subseteq \T\mathcal{X}_i$
  has a \Grobner basis in $\T\mathcal{X}_i$ given by the $f_{i,j}$'s and
  the $g_{a,b}=ab-ba$ for $a>_i b \in \mathcal{X}_i$, it is enough to show
  that for all monomials $m$ and $n$ of $\T\mathcal{X}$, we may write
  $mn-nm-\tilde{\pi}(m,n)=\sum_{r} \mu_r g_{a_r,b_r} \eta_r$ for some $\mu_r,\eta_r \in
  \T\mathcal{X}$  and some $a_r \gg b_r \in \mathcal{X}$.  The statement
  then follows from repeated application of the bi-derivation law.  Indeed, 
  writing $m=m'm''$, we have that:
  \begin{eqnarray*}
    \lefteqn{m'm''n-nm'm''-\tilde{\pi}(m'm'',n)} \\
    &=& m'm''n-nm'm''-m'\tilde{\pi}(m'',n)-\tilde{\pi}(m',n)m''\\
    &=& m'(m''n-\tilde{\pi}(m'',n))-(nm'+\tilde{\pi}(m',n))m'' \\
    &=& m'(m''n-\tilde{\pi}(m'',n)) -m'nm'' + m'nm''-(nm'+\tilde{\pi}(m',n))m'' \\
    &=& m'(m''n-nm''-\tilde{\pi}(m'',n)) + (m'n-nm'-\tilde{\pi}(m',n))m'' 
  \end{eqnarray*}
  which 
  allows us to continually reduce the length of the monomials under
  consideration, so that we may assume that $m$ and $n$ both have length
  $1$, i.e. that they may be identified as elements of $\mathcal{X}$.
  Then we may write either $mn-nm-\tilde{\pi} = g_{m,n}$ if $m\gg n$ or otherwise
  $mn-nm-\tilde{\pi} = -g_{n,m}$.

  We now show that this is indeed a \Grobner basis.
  Without  loss of generality, we may assume that the leading coefficient 
  $f_{i,j}$  is $\LC(f_{i,j})=1$.   Recalling that the $f_{i,j}$ are
  polynomials in $\mathbf{S}\mathcal{X}_i\rInto \T\mathcal{X}_i$,
  we may write their leading monomials as:
  \begin{displaymath}
    \LM(f_{i,j})=:\uvec{x}_{i,j}=x_{i,j,1}x_{i,j,2}\cdots x_{i,j,n_j}
  \end{displaymath}
  where the $x_{i,j,k} \in \mathcal{X}_i$ satisfy
  $x_{i,j,1} \le_i x_{i,j,2} \le_i \cdots \le_i x_{i,j,n_j}$.
  For the polynomials $g_{y,x}$ with $y>x \in \mathcal{X}$, the strong filtration
  property on $\pi$ ensures that $\LM(g_{y,x})=yx$.
  Now it suffices to check the vanishing of the following
  ``S-polynomials'':
  \begin{displaymath}
    S_1 = g(x_k,x_{i,j,1})x_{i,j,2}\cdots x_{i,j,n_j} - x_k f_{i,j}
  \end{displaymath}
  with $x_k \gg x_{i,j,1}$ in $\mathcal{X}$:
  \begin{displaymath}
    S_2 = x_{i,j,1}x_{i,j,2}\cdots x_{i,j,n_j-1}g(x_{i,j,n_j},x_k)-  f_{i,j} x_k
  \end{displaymath}
  with $x_{i,j,n_j} \gg x_k$ in $\mathcal{X}$,
  and:
  \begin{displaymath}
    S_3(k,j,i) = g(x_k,x_j)x_i - x_k g(x_j,x_i)
  \end{displaymath}
  with $x_k\gg x_j \gg x_i$ in $\mathcal{X}$.

  Let us write $f_{i,j} = \uvec{x}_{i,j} - h_{i,j}$.  
  We have two possibilities in showing that $S_1$ vanishes.
  The first case is for $x_k \in \mathcal{X}_i$, in which we have for
  all $1 \le r \le n_j$:
  \begin{displaymath}
    g(x_k,x_{1,j,r})= x_kx_{i,j,r}-x_{i,j,r}x_k
  \end{displaymath}
  by internal symmetry.  From this, the vanishing of $S_1$ follows easily.
  The second case is for $x_k \not\in\mathcal{X}_i$, but as $x_k \gg
  x_{i,j,1}$,  this  means that $x_k \in \mathcal{X}_{i'}$ for some $i'> i$,
  so that $x_k \gg x_{i,j,1} \gge
  g_{i,j_2}\cdots \gge  x_{i,j,n_j}$.
  We also note that as
  $\tilde{\pi}$ descends to $\pi$ on $A$ we necessarily have that
  for all $h \in \T\mathcal{X}$ we have that:
  \begin{displaymath}
    0=\tilde{\pi}(f_{i,j},h) = \tilde{\pi}(\uvec{x}_{i,j} - h_{i,j},h).
  \end{displaymath}
  Now as $\tilde{\pi}$ is a bi-derivation, it follows that:
  \begin{displaymath}
    \sum_{1\le r \le n_j} x_{i,j,1}\cdots
    x_{i,j,r-1}\tilde{\pi}(x_{i,j,r},h) x_{i,j,r+1}\cdots x_{i,j,n_j}
    = \tilde{\pi}(h_{i,j},h)
  \end{displaymath}
  With these observations,  we may now compute $S_1$ as follows:
  \begin{eqnarray*}
    \lefteqn{S_1 =
    -x_{i,j,1}x_k x_{i,j,2}\cdots x_{i,j,n_j} -
    \tilde{\pi}(x_k,x_{i,j,1})x_{i,j,2}\cdots x_{i,j,n_j}  + x_k
    h_{i,j}} \\
    &=&
    -x_{i,j,1} x_{i,j,2} x_k x_{i,j,3}\cdots x_{i,j,n_j} -
    x_{i,j,1}\tilde{\pi}(x_k,x_{i,j,2})x_{i,j,3}\cdots x_{i,j,n_j}  \\
    &&
    -\tilde{\pi}(x_k,x_{i,j,1})x_{i,j,2}\cdots x_{i,j,n_j}  + x_k
    h_{i,j} \\
    &\vdots & \\
    &=& -x_{i,j,1}\cdots x_{i,j,n_j} x_k     + x_k h_{i,j} \\
    & & - \sum_{1\le r \le n_j} x_{i,j,1}\cdots 
    x_{i,j,r-1},\tilde{\pi}(x_k,x_{i,j,r}) x_{i,j,r+1}\cdots x_{i,j,n_j}
    \\
    &=&
    -h_{i,j}x_k - \tilde{\pi}(x_k,\uvec{x}_{i,j}) + x_k h_{i,j} \\
    &=& \tilde{\pi}(x_k,h_{i,j}) -
    \tilde{\pi}(x_k,\uvec{x}_{i,j}) =0
  \end{eqnarray*}
  The vanishing of $S_2$ is similar.

  The vanishing of $S_3$ proceeds similarly to the case of the \Grobner basis
  for the PBW theorem of a Lie algebra.
  Indeed, we have:
  \begin{eqnarray*}
    S_3 &=& -x_kx_ix_j - x_k\tilde{\pi}(x_j,x_i) + x_jx_kx_i + \tilde{\pi}(x,x_j)x_i \\
    &=& -x_ix_kx_j -\tilde{\pi}(x_k,x_i)x_j-x_k\tilde{\pi}(x_j,x_i) +
    x_jx_ix_k \\
    & & + x_j\tilde{\pi}(x_k,x_i) + (x_k,x_j)x_i \\
    &=&
    -x_ix_jx_k-x_i\tilde{\pi}(x_k,x_j)-\tilde{\pi}(x_k,x_i)x_j-x_k(x_k,x_i)
    + x_ix_jx_k\\
    && +\tilde{\pi}(x_k,x_i)x_k + x_j\tilde{\pi}(x_k,x_i) + \tilde{\pi}(x_k,x_j)x_i \\
    &=& \tilde{\pi}(\tilde{\pi}(x,x_j),x_i) + \tilde{\pi}(x_j,\tilde{\pi}(x_k,x_i)) + \tilde{\pi}(\tilde{\pi}(x_j,x_i),x_k) \\
    &=& 0
  \end{eqnarray*}
  where the last equality follows from skew-Symmetry along with the Jacobi
  identity.
\end{proof}

\begin{eRemark}
\item
  We have injective maps $A_i\rInto^{\iota_{i}} T=\coprod_i A_i$, and it is
  clear   that the composite maps $A_i\rTo^{q\comp \iota_i} U$ is
  injective, so that $A_i$ is a $k$-sub-algebra of $U$.
\item
  Let $\g$ be a Lie algebra  over $k$ and choose a basis
  $\g=k\langle X_i\rangle_{i\in I}$.  Take for $A_i$ the
  free commutative algebra generated by $x_i$.
  Then $T=\T\g$, and $U=\T\g/\mathcal{J}=\U\g$, the universal enveloping
  algebra, and we may identify $S=\bigotimes A_i = \mathbf{S}\g$.
  Choosing an ordering $<$ on $I$,  we recover the usual PBW
  theorem.
\item
  More generally, let $\g$ be a Lie algebra where we may write
  $\g=\oplus_{i\in I} \h_i$ where the $\h_i$ are abelian Lie sub-algebras.
  Taking $A_i = \mathbf{S}\h_i=\U\h_i$, then we have:
  \begin{displaymath}
    T = \T\g /\mathcal{R}.
  \end{displaymath}
  Here, $\mathcal{R}$ is the ideal $\mathcal{R}=\sum_i\mathcal{R}_i$
  where $\mathcal{R}_i$ is the ideal 
  generated by the   elements $a a' - a' a$ for $a,a' \in A_i$.
  Then, the  Lie algebra structure on $\g$ induces a Poisson bracket on $T$, and
  as the $\h_i$ are abelian sub-algebras, we see that the Poisson bracket has
  internal symmetry with respect to the $A_i$.  We may identify
  $U=\U\g$, and we may view the $A_i=\mathbf{S}\h_i$ as sub-algebras of $U$.
  Identifying $S=\bigotimes A_i = \mathbf{S}\g$, and choosing an ordering
  on the indexing set $I$ as well as an ordered basis for each of the
  $\h_i$'s,  we again recover the usual PBW theorem.
\end{eRemark}

\begin{Definition}\label{strong_mult}
  Let $\overline{T}$ be a splay with coordinate algebras $A_i$ for $i\in I$
  and suppose that $\overline{T}$ has a strongly filtered Poisson structure $\pi$.
  Denote the multiplication of $S:=\otimes A_i$ be $\cdot$.  Let 
  $<$ be an ordering on $I$.  Then we say that a map  
  $T\rTo^{\Delta}T \tensor T$ is {\bf strongly multiplicative}
  (with respect to $<$) if:
  \begin{displaymath}
    \Delta \pi = \left( \cdot \tensor \pi + \pi\tensor \cdot + \pi\tensor
    \pi\right )\comp (\1\tensor \tau \tensor \1) \comp \Delta.
  \end{displaymath}
  Here, $\1$ is the identity map, and $\cdot$ denotes, by abuse of notation,
  the map induced by on $\T\mathcal{X}$ by:
  \begin{displaymath}
    \mathcal{X}\times\mathcal{X} \rTo^{\cdot} \mathbf{S}\mathcal{X} \rInto^{\iota_{\ll}} T\mathcal{X}.
  \end{displaymath}
\end{Definition}

\begin{eRemark}
\item
  Suppose that we are in the situation of \Ref{Def}{strong_mult}.  Then if $\eta,\zeta\in\mathcal{X}$
  we have $\eta\cdot\zeta = \zeta\cdot\eta$.  Moreover, if
  $\eta>\zeta$ we have by construction of $T$ that:
  \begin{displaymath}
    \zeta\eta = \zeta\cdot\eta = \eta\cdot eta
  \end{displaymath}
\item
  Let us continue the example given by a Lie algebra $\g$ with the choice
  of an ordered basis, where $\T\g\rTo^{\Delta}\T\g\tensor \T\g$
  is map induced by the diagonal map on $\g$.
  Then, one has that the $\pi\tensor\pi$ term of the
  strong multiplicativity condition automatically vanishes when applied
  to elements of $\g \subseteq T=\T\g$ and, thus, what remains is the
  usual multiplicativity condition (c.f. \cite{LuWei}).
  Moreover, one can view the multiplicativity condition as exactly the
  condition  that is needed to ensure that $\Delta$ descends to an algebra map
  $\U\g\rTo^{\Delta} \U\g\tensor\U\g$.
\end{eRemark}

\noindent In view of the previous remark, we note the following:

\begin{Proposition}
  Let $\overline{T}$ be a splay with Poisson structure $\pi$ and
  coordinate algebras $A_i$ for $i\in I$ satisfy the conditions of \Ref{Lem}{Tbasis}.
  Suppose 
  that for each $i$ we have algebra maps $A_i \rTo^{\Delta_i} A_i\tensor
  A_i$, such that the induced map $T \rTo^{\Delta} T\tensor T$
  is strongly multiplicative.  Then $\Delta$ descends to an algebra
  map $U\rTo^{\Delta} U\tensor U$, where $U$ is as in \Ref{Thm}{PBW}.
\end{Proposition}
\begin{proof}
  Let us fix an ordering on the index set $I$.
  In the proof of  \Ref{Thm}{PBW} we realized $U$ as the quotient of
  $\T\mathcal{X}$ by $\mathcal{R}+\tilde{\mathcal{J}}_{\ll}$ for which we
  have exhibited  a \Grobner
  basis for $\mathcal{R} + \tilde{\mathcal{J}}_{\ll}$ given in terms of the $f_{i,j}$'s
  and  $g_{y,x}$'s for $y\gg x \in \mathcal{X}$.  Let us write $\pi$ instead
  of $\tilde{\pi}$ by abuse of notation.
  Since $\Delta$ is induced from the $\Delta_i$,
  it is enough to show that  $\Delta$ vanishes on the $g_{y,x}$'s.
  Let us suppose that $y\in \mathcal{X}_j$ and $x\in\mathcal{X}_i$.
  If $i=j$, then $g_{y,x}=yx-xy$ and the assertion is trivial.
  If $i\ne j$, then $j > i$ as $y\gg x$.
  As $\Delta$ was induced from the $\Delta_i$'s, we may write:
  \begin{eqnarray*}
    \Delta(x) &=& \Delta_i(x) = \sum_r \lambda_r x'_r\tensor x''_r 
    \\
    \Delta(y) &=& \Delta_j(y) = \sum_s \mu_s y'_s\tensor y''_s.
  \end{eqnarray*}
  Here, the     $x'_r,x''_r$ are  monomials in $\mathcal{X}_i$,
  and the $y'_s,y''_s$ are  monomials in $\mathcal{X}_j$.
  We compute modulo the ideal
  $(\mathcal{R}+\tilde{\mathcal{J}}_{\ll})\tensor \T\mathcal{X} + \T\mathcal{X}\tensor
  (\mathcal{R}+\tilde{\mathcal{J}}_{\ll})$ in 
  that:
  \begin{eqnarray*}
    \lefteqn{\Delta(yx-xy) \equiv
      \sum_{r,s} \lambda_r \mu_s
      \left(
        y'_sx'_r \tensor y''_s x''_r -
        x'_ry'_s\tensor x''_ry''_s 
      \right)
    }\\
    &\equiv&\sum_{r,s} \lambda_r\mu_s
    \left(
      (x'_ry'_s +\pi(y'_s,x'_t))\tensor(x''_ry''_s + \pi(y''_s,x''_r))-
      x'_ry'_s\tensor x''_ry''_s 
    \right)\\
    &\equiv&
    \sum_{r,s} \lambda_r\mu_s
    \left(\begin{array}{l}
      \pi(y'_s,x'_r)\tensor x''_ry''_s  + x'_ry'_s\tensor \pi(x''_r,y''_s)\\
      + \pi(y'_s,x'_r)\tensor \pi(y''_s,x''_r)
      \end{array}
    \right)
    \\
    &\equiv&
    (\pi\tensor \cdot + \cdot\tensor \pi + \pi\tensor \pi) \comp (1\tensor
    \tau \tensor 1) \comp (\Delta\tensor \Delta) (y\tensor x)
    \\
    &\equiv& \Delta \comp \pi (y,x) 
  \end{eqnarray*}
  which yields the desired result.
\end{proof}

\section{Formal Groups and Distribution Algebras}

We wish to  apply the results of $\S$\ref{PoissonAlgebras} in order to give a description of
non-commutative formal groups over a field $k$ of characteristic $p > 0$.
We  establish some notation and conventions for general
rings $R$. We say that $\GG$ is a formal group (over $R$) if it is a group object
in the category of smooth formal varieties over $R$.  We denote its dimension
by $n<\infty$.

We recall some basic results of such formal groups (c.f \cite{DemGab}, \cite{Dieudonne},
\cite{Hels}, \cite{Jan}, and \cite{Lazard}).  
Given a  formal group $\GG$ over $R$, one can look at
the ring of formal functions $R[\![\GG]\!]$ whose maximal ideal of functions
we denote $\mathfrak{m}$.   One can form the  continuous linear dual of $R[\![\GG]\!]$
with respect to the $\mathfrak{m}$-adic topology which we denote by
$\Dist(\GG)$.  We denote by $\langle \ \ , \ \ \rangle$ the natural
pairing between $\Dist(\GG)$ and $R[\![\GG]\!]$.
One
knows that $\Dist(\GG)$ has an algebra structure
given by dualizing the co-multiplication map $m$ of
$R[\![\GG]\!]$.  The non-commutativity of $\GG$ is reflected equally in the
non-commutativity of the algebra structure of $\Dist(\GG)$ and
the non-co-commutativity of the co-multiplication $m$.   One knows further that
one may dualize the multiplication of $R[\![\GG]\!]$ to give
a co-multiplication map $\Delta$ on $\Dist(\GG)$.  The algebra and co-algebra
structure of $\Dist(\GG)$, and also those of $R[\![\GG]\!]$, are compatible
so that we in fact have a bi-algebra structure.  
The $\mathfrak{m}$-adic topology on $R[\![\GG]\!]$ yields a filtration of $\Dist(\GG)$
by sub-co-algebras:
\begin{eqnarray*}
  \Dist(\GG) &=& \bigcup_{n\ge 0 } \Dist^{(n)}(\GG) \\
  \Dist^{(n)}(\GG) &:=& R\oplus \Dist^{(n)}_{+}(\GG) \\
  \Dist^{(n)}_+(\GG) &:=& \Hom_R(\mathfrak{m}/\mathfrak{m}^{n+1},R).
\end{eqnarray*}
A choice of coordinate system $\uvec{x} = \{x_1,\cdots,x_n\}$ on $\GG$
gives an identification $R[\![\GG]\!] = R[\![x_1,\cdots,x_n]\!]$.  If one
chooses an ordering on the $x_i$'s, say $x_1<\cdots <x_n$, then
then one is  working in the category of formal group laws.
Any   choice of a coordinate system,  without any specified ordering, gives
the choice of an additional structure on $\Dist(\GG)$, namely that of 
a {\bf DVPS-bi-algebra}.  By this we mean there is a choice of a basis
$\delta_{\uvec{x}^I}$ of $\Dist(\GG)$ where $J=(J_1,\cdots, J_n)\in \ZZ_{\ge 0}^n$
is a multi-index such that  the co-multiplication law satisfies:
\begin{displaymath}
  \Delta(\delta_{\uvec{x}^J}) = \sum_{A+B=J} \delta_{\uvec{x}^A} \tensor \delta_{\uvec{x}^B}.
\end{displaymath}
The specific basis that the choice of $\uvec{x}$ is the basis given
by the maps:
\begin{displaymath}
  \delta_{\uvec{x}^J} (\uvec{x}^K) := \left\{ \begin{array}{cc}1 & J=K \\ 0 & J \ne K \end{array}\right.
\end{displaymath}
where for the multi-index $K$ we have defined  $\uvec{x}^K:=x_1^{K_1}\cdots x_n^{K_n}
\in R[\![\GG]\!]$.  We will call such a basis an {\bf additive basis}.
The existence of the additive basis means that one can recognize the dual of
$\Dist(\GG)$ with multiplication given by dualizing $\Delta$ as a ring of
power series.
There is a final structure on $\Dist(\GG)$ which is the
anti-pode or inverse $\Dist(\GG)\rTo^{\inv} \Dist(\GG)$
which gives $\Dist(\GG)$ the structure of a Hopf
algebra.  The existence and uniqueness is automatic in the situations we
will consider below.

\begin{eExample}
\item\label{example_of_g}
  Let  $\mathbf{G}$ be a smooth (analytic, algebraic,..)  group defined over $R$.
  Then we denote the completion of its ring of functions at the identity by
  $R[\![\mathbf{G}]\!]$,   which we  identify 
  as the ring of formal functions  of a formal group.  Denote by $\g$ the 
  Lie algebra of $\mathbf{G}$.
  If $R$ is a $\QQ$-algebra, then $\U\g$ is
  canonically isomorphic to $\Dist(\mathbf{G})$.
\item
  If $k$ is a field of characteristic $p >0$, or more generally $R$ is an
  $\FF_p$-algebra, then 
  one has that $\Dist(\GG)$ is filtered by sub-bi-algebras:
  \begin{displaymath}
    \Dist(\GG) = \bigcup_{r\ge 0} \Dist^{(p^r -1)}(\GG).
  \end{displaymath}
\end{eExample}

\noindent Let us now assume that $R=k$ is a field of characteristic $p>0$.
In order to make use of the results of $\S1$, we will
need to make use of a second basis:

\begin{LemmaDefinition}
  Let $\GG$ be a formal group over a field $k$ of characteristic $p>0$.  Let
  $\uvec{x}_< = \{x_1<\cdots < x_n\}$ be the choice of an ordered coordinate
  system for $\GG$.  Then $\Dist(\GG)$ has a $k$-basis, the
  {\bf multiplicative basis},  given  for multi-indices $J$
  as  the ordered products:
  \begin{displaymath}
    {\delta}_{\uvec{x}_<}^{\uvec{J}} = {\delta}_{x_1}^{\uvec{J_1}}
    \cdots {\delta}_{x_n}^{\uvec{J_n}} 
  \end{displaymath}
  where:
  \begin{displaymath}
    {\delta}_{x_j}^{\uvec{J_j}} :=
    \delta_{x_j^{p^0}}^{J_{j,0}}\cdots\delta_{x_j^{p^r}}^{J_{j,r}} .
  \end{displaymath}
  Here, the $\delta_{x_j^{p^s}}$'s are as in the additive basis above, 
  and the $J_{j,t}$'s are given by the $p$-adic expansions of the $J_j$'s:
  \begin{displaymath}
    J_j = J_{j,0} p^0 + J_{j,1} p^1 + \cdots + J_{j,r} p^r.
  \end{displaymath}
  When the choice of an ordered coordinate system $\uvec{x}_<$ is
  understood,   we will drop the subscripts and simply write
  ${\delta}^{\uvec{J}}:={\delta}^J_{\uvec{x}_<}$.
\end{LemmaDefinition}

\begin{eRemark}
\item\label{multiplicativebasis}
  The fact that this is indeed a basis is in  \cite{Dieudonne}.
  Defining:
  \begin{displaymath}
    \mathcal{X}:= \{ \delta_{x_j}^{p^r} \ | \ 1\le j \le n \ \ r\ge 0 \} 
  \end{displaymath}
  one has at once an 
  induced $R$-linear map $\iota_{\uvec{x}_<}$  given by:
  \begin{displaymath}
    \begin{diagram}[height=1.7em]
      \Dist(\GG) & \rInto^{\iota_{\uvec{x}_<}} & \T\mathcal{X} \\
      \delta^{\uvec{I}} &\rMapsto       \delta^{\uvec{I}}.
    \end{diagram}
  \end{displaymath}
  We define for $\eta,\zeta \in \mathcal{X}\subseteq \Dist(\GG)$:
  \begin{eqnarray*}
    \tilde{\pi}_{\uvec{x}_<}(\eta,\zeta)&:=& \iota_{\uvec{x}_<}(\eta\zeta -
    \zeta\eta) =\iota_{\uvec{x}_<} \comp \pi_{\canon}(\eta,\zeta)\\
    g_{\uvec{x}_<,\eta,\zeta}&:=& \iota_{\uvec{x}_<}(\eta)\iota_{\uvec{x}_<}(\zeta) -
    \iota_{\uvec{x}_<}(\zeta)\iota_{\uvec{x}_<}(\eta) -
    \tilde{\pi}_{\uvec{x}_<}(\eta,\zeta) \\
    F_{\uvec{x}_<}(\eta) &:=& \iota_{\uvec{x}_<}(\eta^p) \\
    f_{\uvec{x}_<,\eta} &:=& \iota_{\uvec{x}_<}(\eta)^p - F_{\uvec{x}_<}(\eta).
  \end{eqnarray*}
  The extension of $\tilde{\pi}_{\uvec{x}_<}$ by the bi-derivation law to $\T\mathcal{X}$
  will also be denoted by $\tilde{\pi}_{\uvec{x}_<}$.
  If we take the choice $\uvec{x}_<$ as understood
  we shall
  write $\tilde{\pi} = \tilde{\pi}_{\uvec{x}_<}$.  We may also
  identify
  $\mathcal{X}\subseteq\Dist(\GG)$ and $\Dist(\GG)$
  with their images in $\T\mathcal{X}$ via $\iota_{\uvec{x}_<}$.  Thus
  we may write more succinctly:
  \begin{eqnarray*}
    g_{\eta,\zeta} &=& \eta\zeta -\zeta\eta - \tilde{\pi}(\eta,\zeta) \\
    f_{\eta} &=& \eta^p - F(\eta) 
  \end{eqnarray*}
\item\label{properties_of_relations}
  We have the following important filtration properties of the $F(\eta)$ and
  $\pi_{\canon}(\eta,\zeta)$:
  \begin{eqnarray*}
    F(\delta_{x_i^{p^r}}) &\in&  \Dist^{(p^{r+1}-1)}(\GG) \\
    \pi(\delta_{x_j^{p^s}},\delta_{x_i^{p^r}}) &\in& \Dist^{(p^r+p^s -1)}(\GG).
  \end{eqnarray*}
\item
  There is an analogue to the multiplicative basis in the characteristic zero situation of $\g$
  as above.  For simplicity, let us assume that $\g$ is abelian and arises
  as the Lie algebra of the formal group $\GG=\GG_a^n$.
  Let us choose a coordinate
  $y_i$ for each of the  $\GG_a$'s  such that  we may write the co-multiplication
  map  of $\GG$ as:
  \begin{displaymath}
    m(y_i) = y_i \tensor 1 + 1 \tensor y_i.
  \end{displaymath}
  Let us denote the corresponding basis of $\g$ by $g=\langle Y_i
  \rangle_{1\le i \le n}$
  where $Y_i(y_j)=\delta_{i,j}$ is the Kronecker delta product.
  Then, in this setting the, ``multiplicative
  basis'' of $\U\g=\mathbf{S}\g=\Dist(\GG)$ is given by the elements:
  \begin{displaymath}
    Y^{\uvec{I}} = Y_1^{I_1}\cdots Y_n^{I_n}.
  \end{displaymath}
  We
  may compute:
  \begin{eqnarray*}
    \delta_{y_i} \delta_{y_i^k}  &=& \sum_{J} \langle m(\uvec{y}^J),
    \delta_{y_i}\tensor \delta_{y_i}^k \rangle \delta_{\uvec{y}^J} \\
    &=& \sum_{j} \langle (y_i\tensor 1 + 1\tensor y_i)^j,
    \delta_{y_i}\tensor \delta_{y_i}^k \rangle \delta_{y_i^j} \\
    &=& (k+1) \delta_{y_i^{k+1}}
  \end{eqnarray*}
  so that  identifying $Y_i = \delta_{y_i}$
  the additive and multiplicative basis are related by:
  \begin{displaymath}
    {Y}^{\uvec{I}} = \delta_{\uvec{y}}^{\uvec{I}} = {I!}\delta_{\uvec{y}^{I}}
  \end{displaymath}
  where $I!:=I_1!\cdots I_n!$.

  We may have also chosen the coordinates $x_i = \exp(y_i) - 1$ for $\GG$,
  so  that:
  \begin{displaymath}
    m(x_i) = x_i\tensor 1 + 1 \tensor x_i + x_i \tensor x_i
  \end{displaymath}
  which is the group law for $\mathbf{G}_m$ the multiplicative group as well as its formal
  completion $\GG_m$.
  Then as we have that:
  \begin{displaymath}
    \delta_{x_i}\cdot \delta_{x_i^k} = (k+1)\delta_{x_i^{k+1}} + k \delta_{x_i^{k}}
  \end{displaymath}
  one sees that the change of basis between the multiplicative and additive
  bases can be expressed in terms of the Sterling numbers of the first kind:
  \begin{displaymath}
    \delta_{x_i^n} = \frac{1}{n!}\delta_{x_i}(\delta_{x_i} -1)(\delta_{x_i}-2)\cdots(\delta_{x_i} -(n-1)).
  \end{displaymath}
\item
  The impossibility of inverting $I!$ and the fact that $\GG_a$ and $\GG_m$
  are  not isomorphic in   positive   characteristic
  is the essential source of all the complications.
  Let us define for $n\in \ZZ\ge 0$:
  \begin{displaymath}
    n!_p := n_0!n_1! \cdots n_r!
  \end{displaymath}
  where $n=n_0 + n_1 p + \cdots n_r p^r$ is the $p$-adic expansion of $n$.
  We make a similar definition for a multi-index $I$.
  Then  one know that for a formal group $\GG$ we have:
  \begin{displaymath}
    \delta_{\uvec{x}^I} = \frac{1}{I!_p} \delta_{\uvec{x}_<}^{\uvec{I}} + \ \ \mbox{ Lower Order Terms }
  \end{displaymath}
  In the case of $\GG=\GG_a$ with coordinate $y$ as above,  we have, in fact, that:
  \begin{displaymath}
    \delta_{y^n} = \frac{1}{n!_p} \delta_y^{\uvec{n}}.
  \end{displaymath}
  For $\GG=\GG_m$, with coordinate $x$ as above, we have ``Sterling numbers modulo p''.
  Indeed, we  may compute:
  \begin{displaymath}
    \delta_{x^{p^r}}\delta_{x^m} = (m_r +1)\delta_{x^{m+1}} + m_r \delta_{x^m}
  \end{displaymath}
  where $m=m_0 + \cdots + m_r p^r$ is the $p$-adic expansion of $m$ and we assume $m_r < p-1$.
  From this it follows that we have:
  \begin{displaymath}
    \delta_{x^m} = \frac{1}{m!_p} \prod_{t=0}^r \delta_{x^{p^t}}(\delta_{x^{p^t}} -1)(\delta_{x^{p^t}}-2)\cdots(\delta_{x^{p^t}} -(m_t-1)).
  \end{displaymath}
\end{eRemark}

\begin{Proposition}\label{dist_as_quotient}
  Let  $\GG$ be a formal group over a field $k$
  of characteristic $p>0$. Let
  $\uvec{x}_<$ be the choice of an ordered coordinate system.
  Let us adopt the notation  as in \Ref{Rmk}{multiplicativebasis}.
  Define $\deg(\delta_{x_i}^{p^r}):=p^r$ and take:
  \begin{displaymath}
    \delta_{x_j^{p^s}} > \delta_{x_i^{p^r}}
  \end{displaymath}
  if $j>i$ or if $j=i$ and $s > r$.
  Let $\mathcal{I}$ be the two-sided ideal of $\T\mathcal{X}$:
  \begin{eqnarray*}
    \mathcal{I}&=& \mathcal{R} + \tilde{\mathcal{J}}_{\ll} \\
    \mathcal{R}&=& \left( f_{\eta} \ | \  \eta \in \mathcal{X} \right) \\
    \tilde{\mathcal{J}}_{\ll} &=& \left( g_{\eta,\zeta} \ | \ \eta > \zeta \in \mathcal{X} \right).
  \end{eqnarray*}
  Imposing the graded lexicographic ordering on $\T\mathcal{X}$, the set:
  \begin{displaymath}
    \{f_{\eta},g_{\eta,\zeta} \ | \eta > \zeta \in \mathcal{X} \}
  \end{displaymath}
  is a \Grobner basis for the ideal $\mathcal{I}$.  In particular, $\Dist(\GG)$
  is isomorphic as an algebra to $\T\mathcal{X}/\mathcal{I}$.
\end{Proposition}

\begin{proof}
  Due to the monomial ordering indicated and by
  \Ref{Rmk}{properties_of_relations}
  we see that we having the following leading monomials:
  \begin{eqnarray*}
    \LM(f_{\eta}) &=& \eta^p \\
    \LM(g_{\eta,\zeta}) &=& \eta\zeta.
  \end{eqnarray*}
  The result follows as otherwise we would not have the existence of a multiplicative basis.
\end{proof}

\section{Geometric Formal Groups}

Although \Ref{Prop}{dist_as_quotient} gives some reasonable control on the
algebra structure of $\Dist(\GG)$,  some difficulties remain in trying to characterize
formal groups in such a manner.  First, our construction
depended on the choice of an ordered coordinate system.  Second,
we have yet to address the DVPS-co-algebra structure and its
compatibility with the algebra structure.   For this, we will need to
make an assumption on $\GG$, namely that it is ``geometric'' which loosely
means that one can choose a coordinate system of commutative formal subgroups.
One may think of this as specifying a ``local  abelian structure''.  In the
case where the ground ring is a field of characteristic zero, 
the fact that any commutative formal group is isomorphic to a product of
$\GG_a$'s means that there is essentially one local abelian  structure.
However,  in the positive characteristic situation, there are
non-isomorphic (commutative) one-dimensional formal groups. Moreover, there
are commutative formal groups
which are not the products of one-dimensional formal groups, even over an
algebraically closed field.

\begin{eDefinition}
\item
  Let $\GG$ be a formal group over $R$.   We say that $\GG$ is
  {\bf geometric}
  if there are commutative formal subgroups $\HH_i$ for $i$ in some indexing
  set $I$
  such that there is an isomorphism as between $\GG$ and the product of
  formal varieties $\prod_i \HH_i$.
  We denote a geometric formal group by the data $\{\HH_i \rInto^{\iota_i}
  \GG\}_{i\in I}$, or more briefly $\{\HH_i,\GG\}_{i\in I}$.
\item
  A {\bf morphism} of geometric formal groups:
  \begin{displaymath}
    \{\HH_i,\GG\}_{i\in I} \rTo^{\Phi} \{\HH'_j,\GG'\}_{j\in J} 
  \end{displaymath}
  is a morphism of formal groups $\GG\rTo^{\Phi} \GG'$ such that
  for each $i\in I$ there is some $j_i \in J$ and some morphism of
  commutative formal groups $\HH_i \rTo^{\Phi_{j_i}} \HH'_j$
  such that $\Phi$ is the  morphism induced by the
  UMP of products of formal varieties as in the following diagram:
  \begin{displaymath}
    \begin{diagram}[height=1.7em]
      \HH_i & \rInto^{\iota_i} &  \GG \\
      \dTo^{\Phi_{j_i}} &  \ldDashto(2,4)^{\exists !} & \\
      \HH'_{j_i} & & \\
      \dTo^{\iota_{j_i}} & & \\
      \GG'
    \end{diagram}
  \end{displaymath}
\item
  Let $\{\HH_i,\GG\}_{i\in I}$ be a
  geometric formal group over $R$.
  We define the
  {\bf underlying commutative geometric group}  as:
  \begin{displaymath}
    \Comm(\HH_i,\GG) := \{\HH_i,\HH\}
  \end{displaymath}
  where $\HH:= \prod_i \HH_i$ is the {\bf underlying commutative group} and
  the product is taken in the category of commutative formal groups.
\item
  Let $\{\HH_i,\GG\}_{i\in I}$ be a
  geometric formal group over $R$.  Let $d_i := \dim \HH_i$ and denote
  by abuse of notation $I=|I|$, so that we may identify $I=\{1,\cdots |I|\}$.
  Suppose that for each $i\in I$, we have chosen an ordered coordinate
  system $\uvec{x}_{i,<_i} = \{x_{i,1} <_i \cdots <_i x_{i,d_i} \}$.  Let us suppose
  moreover that we have chosen an ordering $<$ on $I$.  We will call the
  resulting ordered coordinate system:
  \begin{eqnarray*}
    \uvec{x}_{\ll} &:=& \{x_{1,1} \ll \cdots \ll x_{1,d_1} \ll x_{2,1} \ll \cdots x_{I,1}
    \ll \cdots \ll x_{I,d_I} \}\\
    &=:& \{ x_1 \ll \cdots \ll x_n \}
  \end{eqnarray*}
  a {\bf geometric  coordinate system}.
  We will call the data:
  \begin{displaymath}
    \{\uvec{x}_{\ll},\HH_i,\GG\}_{i\in I}
  \end{displaymath}
  a {\bf geometric formal group law}.
\end{eDefinition}

\begin{eRemark}\label{comultiplication}
\item
  Suppose that  $\{\HH_i,\GG\}$ is a geometric formal group.  Then, 
  each of the $\Dist(\HH_i)$ carries its own co-multiplication map, say
  $\Delta_i$, as well as $\HH$ with its co-multiplication $\Delta$.  However,
  as we have assumed an isomorphism of formal varieties $\GG=\prod_i \HH_i$,
  we have an isomorphism of co-algebras $\Dist(\GG_i)=\prod_i \Dist(\HH_i)$
  under which we have $\Delta = \prod_i \Delta_i$.
\item
  One sees that the $\Dist(\HH_i)$ are sub-bi-algebras of $\Dist(\GG)$ which
  generate $\Dist(\GG)$ as an algebra.  The choice of the  commutative formal
  subgroups $\HH_i$ are not unique and, in particular, in positive
  characteristic may not be isomorphic.
\end{eRemark}

\begin{eExample}
\item
  If $\GG$ is a commutative formal group, then it is geometric with $I=\{1\}$
  and $\HH_1 = \GG$.
\item
  Suppose that $\mathbf{G}$ is a connected and simply-connected Lie group
  over a field $k$.
  If $\g$ is the Lie algebra of $\mathbf{G}$ as in \Ref{Ex}{example_of_g}, then for any choice of basis
  $k\langle X_i \rangle$ we can exhibit the formal completion of $\mathbf{G}$ as a
  geometric formal group via the formal completions of 
  the Lie subgroups $\mathbf{H}_i=\mathbf{exp}(tX_i)$.  We have that the $\Dist(\mathbf{H}_i)$ are
  exactly the sub-bi-algebras $\mathbf{S}\{X_i\}=\Dist(\mathbf{H}_i)$ of $\U\g=\Dist(\mathbf{G})$.
  In particular, in characteristic zero all formal groups are geometric.

\item
  If $k$ is an algebraically closed field of positive characteristic $p$, then
  the formal completion of any reductive algebraic group is geometric
  with formal subgroups either $\GG_a$ or $\GG_m$.
\item
  In general, we may have more than one type of a multiplicative coordinate.
  For example, consider in the case of $\mathbb{SL}_2$, we have coordinates $(w,x,y)$
  of type $(\GG_a,\GG_m,\GG_a)$ given by:
  \begin{displaymath}
    \left( \begin{array}{cc} 1+x & y \\ w & 1+z   \end{array}\right)
  \end{displaymath}
  where $1+z=\frac{1+yw}{1+x}$.
  We also have coordinates $(s,t,u)$ of type $(\GG_a,\GG_a,\GG_a)$ given by:
  \begin{displaymath}
    \left(\begin{array}{cc} 1 & 0 \\ s & 1 
    \end{array}\right)
    \left(\begin{array}{cc} 1 +t & t \\ -t & 1-t 
    \end{array}\right)
    \left(\begin{array}{cc} 1 & u \\ 0 & 1 
    \end{array}\right).
  \end{displaymath}
\end{eExample}

\begin{Remark}\label{geom_relations}
  Let $R=k$ be a field of positive characteristic.
  Let $\{\uvec{x}_{\ll},\HH_i,\GG\}_{i\in I}$ be a geometric formal group law.
  We may apply \Ref{Prop}{dist_as_quotient} to each of the commutative
  formal subgroups $\HH_i$.  If we
  define for $i \in I$:
  \begin{displaymath}
    \mathcal{X}_i := \{ \delta_{x_{i,k}^{p^r}} \ | \ 1\le k \le d_i \ \ r
    \ge 0 \}
  \end{displaymath}
  then $\mathcal{X} = \coprod_i \mathcal{X}_i$. Fixing an ordering $<$ on
  $I$, let us write
  $\iota_{\uvec{x}_{\ll}}:=\iota_{\ll}$ and
  $\tilde{\pi}_{\uvec{x}_<}:=\tilde{\pi}_{\ll}$
  if there is no danger of confusion.
  As the $\HH_i$'s are commutative formal subgroups, we have that:
  \begin{displaymath}
    F(\mathcal{X}_i) \subseteq \Dist(\HH_i)
  \end{displaymath}
  in other words:
  \begin{displaymath}
    f_{\eta} \subseteq \T\mathcal{X}_i \subseteq \T\mathcal{X} \ \ \ \  \forall
    \eta \in \mathcal{X}_i.
  \end{displaymath}
  Further, we have that for $\eta,\zeta \in \mathcal{X}_i$ that:
  \begin{displaymath}
    g_{\eta,\zeta} = \eta\zeta - \zeta\eta \in \T\mathcal{X}_i \subseteq \T\mathcal{X}
  \end{displaymath}
  \myqed
\end{Remark}

\begin{Lemma}\label{pi_is_sfsm}
  Let us adopt the notation and assumption of \Ref{Rmk}{geom_relations}.
  Then $\tilde{\pi}_{\ll}$ descends to
  a strongly filtered, strongly
  multiplicative Poisson structure $\pi_{\ll}$ on the splay $\overline{T}$
  of the $\Dist(\HH_i)$.
\end{Lemma}

\begin{proof}
  By the presentation of $\Dist(\GG)$ as given in
  \Ref{Prop}{dist_as_quotient},
  and since we have $\tilde{\pi}_{\ll} (\mathcal{X}_i,\mathcal{X}_i) = 0$, 
  we see that $\tilde{\pi}_{\ll}$ descends 
  as we may present $T$ as:
  \begin{displaymath}
    T = \T\mathcal{X} / \sum_i \left( f_{\eta}, g_{\eta,\zeta} \ | \
    \eta,\zeta\in\mathcal{X}_i \right).
  \end{displaymath}
  The internal symmetry follows.
  The strong filtration
  property follows by \Ref{Rmk}{properties_of_relations}.
  As $\Dist(\GG)$ is in fact a bi-algebra, one has by
  \Ref{Rmk}{comultiplication} the strong multiplicativity.

\end{proof}

\begin{Proposition}
  Let $k$ be a field of  characteristic $p>0$.
  Let $\{\HH_i,\GG\}_{i\in I}$ be a geometric formal group over $k$.
  Suppose that $\uvec{x}_{\ll}$  (resp. $\uvec{x'}_{\ll'}$) is a geometric coordinate
  systems on $\GG$ induced by the choice of an ordering $<$ (resp. $<'$) on
  $I$ and ordered coordinate systems $\uvec{x}_{i,<_i}$
  (resp. $\uvec{x}'_{i',<'_i}$)  on the  $\HH_i$.
  Let $\pi$ (resp. $\pi'$)  be the resulting  Poisson structure on $T$ as in
  \Ref{Lem}{pi_is_sfsm}.   Then $\pi$ and $\pi'$ are equivalent.
  In particular if $<=<'$,  then $\Phi=\1$.

  \noindent If 
  $\{\HH_i,\GG\}_{i \in I} \rTo^{\Phi}\{\HH_j',\GG'\}_{j \in J}$
  is a morphism of geometric
  formal groups, then there is an induced Poisson morphism
  $\overline{T}\rTo^{\Phi} \overline{T}'$.
\end{Proposition}

\begin{proof}
  Let us first suppose that we are working with a fixed ordering $<=<'$ on $I$.
  Denote by $\Phi_i$ the isomorphisms:
  \begin{displaymath}
    \Dist(\HH_i)= k\langle \delta_{\uvec{x}_i}^{\uvec{J_i}} \rangle
    \rTo^{\Phi_i} k\langle \delta_{\uvec{x}'_i}^{\uvec{J'_i}} \rangle=\Dist(\HH_i)
  \end{displaymath}
  given by the choices of the two multiplicative bases for the $\Dist(\HH_i)$.
  Adopting the notation of \Ref{Lem}{dist_as_quotient} for each of the coordinate systems, we
  denote the respective generating sets by $\mathcal{X}=\coprod
  \mathcal{X}_i$ and
  $\mathcal{X}'=\coprod \mathcal{X}'_i$.
  The choice of the different coordinate systems gives two different
  presentations   of  $T$ as:
  \begin{eqnarray*}
    T  &=& \T\mathcal{X} / \mathcal{R} +
    \sum_i \left(\eta\zeta-\zeta\eta  \ |  \ \eta,\zeta \in \mathcal{X}_i
    \right) \\
    T  &=& \T\mathcal{X}' / \mathcal{R}' +
    \sum_i \left(\eta'\zeta'-\zeta'\eta'  \ |  \ \eta',\zeta' \in \mathcal{X}_i \right)
  \end{eqnarray*}
  The $\Phi_i$ induce an isomorphism between these two
  presentations.
  The Poisson structures under consideration 
  $\pi$ and $\pi'$  are induced respectively from:
  \begin{align*}
    \tilde{\pi}  := \iota \comp \pi_{\canon} & \mbox{ for } &     \Dist(\GG)  \rInto^{\iota}  \T\mathcal{X} \\
    \tilde{\pi}  := \iota' \comp \pi_{\canon} & \mbox{ for } &     \Dist(\GG)  \rInto^{\iota'}  \T\mathcal{X}'.
  \end{align*}
  Our  aim is to show that $\Phi\comp \pi = \pi' \comp (\Phi\tensor
  \Phi)$.
  First we suppose that we have fixed coordinate systems $\uvec{x_i}
  =\uvec{x_i}'$ on each of the $\HH_i$ so that 
  we have only chosen to change the ordering of $<_i$ to $<_i'$.
  As the induced  multiplicative expansions of $\Dist(\GG)$ are equal up to
  recordings within each of the $\mathcal{X}_i$, we see that the induced map $\Phi$ simply
  induces a change of ordering of the monomials in the $\mathcal{X}_i$.
  The result follows in this setting.

  Let us drop the assumption that $\uvec{x}_i = \uvec{x}'_i$.
  By what we have just said, after a possible reordering
  of the coordinates for each of the $\Dist(\HH_i)$, we may assume that
  that the orderings $<_i$ and $<'_i$ are the same.
  Recalling that
  $\pi$ and $\pi'$ are defined in terms of $\pi_{\canon}$,
  where we have:
  \begin{eqnarray*}
    \Phi\comp\pi_{\canon}(\eta,\zeta) &=&
    \Phi(\eta)\Phi(\zeta)-\Phi(\eta)\Phi(\zeta) \\
    &=& \pi_{\canon}(\Phi(\eta),\Phi(\zeta))
  \end{eqnarray*}
  the result now follows by \Ref{Lem}{measure}, by
  noticing that 
  the multiplicative expansion of $\Dist(\GG)$ in terms of the coordinates
  $\uvec{x}'$ can be obtained from multiplicative expansion of $\uvec{x}'$
  by applying $\Phi$ and the relations $\mathcal{R}'$.

  Now we only need to illustrate the effect of a change of ordering $<$ on
  $I$ to $<'$.
  Let us suppose that the ordering  $<$ is the usual ordering $1<2 < \cdots < |I|$.
  Then, we may view any other ordering $<'$ on $I$ as being given by $<'=<^{\sigma}$ for
  some permutation  $\sigma \in \Sigma_I$.  As the permutations are generated by the transpositions,
  we only need consider the case $\sigma = (i \ \ i\!+\!1)$.
  In this case, we consider the automorphism $\Phi$ of $T$ induced by 
  $\Dist(\HH_j)\rTo^{\Phi_j}\Dist(\HH_j)$ where:
  \begin{displaymath}
    \Phi_j =\left\{ \begin{array}{cc} \inv_j & j \in \{i,i+1\} \\ \1 & \mbox{otherwise}. \end{array}\right.
  \end{displaymath}
  Here, $\inv_j$ is the inverse of $\HH_j$.
  The fact that this gives an equivalence of Poisson structures follows from
  the fact that $\inv$ is an anti-algebra homomorphism so
  that for $\eta,\zeta \in \mathcal{X}_i\cup \mathcal{X}_{i+1}$:
  \begin{eqnarray*}
    \inv\comp \pi_{\canon}(\eta,\zeta) &=& \inv(\eta\zeta-\zeta\eta) \\
    &=& \inv(\eta)\inv(\zeta) -\inv(\zeta)\inv(\eta) \\
    &=& \pi_{\canon} (\inv(\zeta),\inv(\eta)).
  \end{eqnarray*}
  To obtain the morphism $\overline{T}\rTo^{\Phi}\overline{T}'$, we
  first choose an ordering $<'$ on the set $J$.  Then we may choose
  any ordering $<$ on $I$ such that if $i_1 < i_2$ then $j_{i_1} <' j_{i_2}$
  where the $j_i$'s are as in the definition of a geometric morphism.
  We note that such an ordering exists:  We fix a choice of such $j_i$'s
  and define $I_j := \{ i\in I \ | \ j_i = j\}$.  Then as we have
  $I =  \coprod_{j\in J}  I_j$, we choose any ordering $<_j$ on the $I_j$
  and define the ordering $<$ on $I$ as $i_1 < i_2$  if either $j_{i_1} < j_{i_2}$
  or $j:=j_{i_1}=j_{i_2}$ and $i_1 <_j i_2$.
  The result follows readily.
\end{proof}

\begin{Theorem}
  Let $\{ \HH_i,\HH\}$ be a commutative geometric formal group.  Suppose
  that $\pi$ is
  a strongly filtered, strongly
  multiplicative Poisson structure on the splay $\overline{T}$ of the
  $\Dist(\HH_i)$.  
  Let $\mathcal{J}$ be the ideal in the splay algebra $T$: 
  \begin{displaymath}
    \mathcal{J} = (fg -gf -\pi(f,g) \ | \ f,g \in T).
  \end{displaymath}
  Then the quotient:
  \begin{displaymath}
    U:=T/\mathcal{J}
  \end{displaymath}
  is the distribution algebra of a formal group, $\GG$. Moreover, 
  $\{\HH_i,\GG\}$ is a geometric formal group.
  
  \noindent Suppose that $\overline{T}'$ is the splay of a
  commutative geometric formal group $\{ \HH'_i,\HH'\}$ with a strongly
  filtered,   strongly multiplicative Poisson structure $\pi'$.  Suppose
  that $\overline{T}\rTo^{\Phi}\overline{T}$ is a Poisson morphism, then
  there is morphism of geometric formal groups
  $\{\HH_i,\GG\} \rTo^{\Phi}\{\HH'_i, \GG'\}$.
\end{Theorem} 

\begin{proof}
  Let us choose an ordering $<$ on $I$.
  By the work that we have done in  $\S1$ we have a bi-algebra
  structure on $U$.   Thus, we only need
  to show that there is a DVPS algebra structure on $U$.  Let us write:
  \begin{displaymath}
    \{ x_1 \ll \cdots \ll x_n\} := \{ x_{1,1} \ll \cdots \ll x_{|I|,d_{|I|}} \}.
  \end{displaymath}
  For a multi-index $J$, we define the following element of $U$:
  \begin{displaymath}
    \delta_{x^J} := \delta_{x_1^{J_1}}\cdots \delta_{x_n^{J_n}}.
  \end{displaymath}
  If $1\le j \le n$, let $i_j \in I$ denote the index such that $x_j$ is a coordinate
  of $\HH_{i_j}$.
  Then, as we have taken care of the ordering,  and by \Ref{Rmk}{comultiplication}, we have:
  \begin{eqnarray*}
    \Delta(\delta_{x^J}) &=& \Delta_{i_1}(\delta_{x_1^{J_1}})\cdots \Delta_{i_n}(\delta_{x_n^{J_n}}) \\
    &=& \left(\sum_{a_1+b_1 =J_1} \delta_{x_1^{a_1}}\tensor \delta_{x_1^{b_1}}\right)
    \cdot \left(\sum_{a_n+b_n =J_n} \delta_{x_n^{a_n}}\tensor \delta_{x_n^{b_n}}\right)\\
    &=& \sum_{a_j+b_j =J_j} \delta_{x_1^{a_1}}\cdots\delta_{x_n^{a_n}}\tensor \delta_{x_1^{b_1}}\cdots\delta_{x_n^{b_n}} \\
    &=& \sum_{A+B = J} \delta_{x^A}\tensor \delta_{x^B}
  \end{eqnarray*}
  as desired.
  The fact that $\Phi$ is a morphism of formal groups is now obvious.
\end{proof}

\begin{Example}
  We wish to take some time to describe the relevant features of the above theory in terms of
  the formal group corresponding to $\mathbb{T}_2$,
  the upper triangular $2\times 2$ matrices of determinant one.   Let us choose a coordinate system
  as follows:
  \begin{displaymath}
    \left( \begin{array}{cc} 1+x & y \\ 0 & (1+x)^{-1} \end{array}\right)
  \end{displaymath}
  by this we mean that the co-multiplication laws, in these coordinates, are given by:
  \begin{eqnarray*}
    x &\rMapsto^{m} & x\tensor 1 + 1 \tensor x + x\tensor x \\
    y &\rMapsto^{m} & (1+x)\tensor y + y\tensor (1+x)^{-1} \\
    && = (1+x)\tensor y + y\tensor (1-x+x^2 - x^3 + \cdots ).
  \end{eqnarray*}
  Then $x$ and $y$ are geometric coordinates for the subgroups $\GG_m$ and $\GG_a$, respectively.
  Although we know that the underlying commutative group is given by $\GG_m \times \GG_a$, the
  PBW-theorem gives in fact a different set of coordinates for this group.  Here, we need
  to consider the two choices of orderings $x<y$ and $x>'y$, and then a calculation shows 
  that the resulting PBW-coordinates are:
  \begin{displaymath}
    \begin{array}{cc}
      \left( \begin{array}{cc} 1+x & y \\ 0 & (1+x) \end{array}\right) &
      \left( \begin{array}{cc} (1+x)^{-1} & y \\ 0 & (1+x)^{-1} \end{array}\right) \\
      x < y  & x >'y. 
    \end{array}
  \end{displaymath}

  We give a complete description of $\mathbb{T}_2$.  
  From the commutative theory, one knows that for $\GG_m$ we have  $\delta_{x^{p^r}}^p = \delta_{x^{p^r}}$
  and for $\GG_a$ we have $\delta_{y^{p^r}}^p = 0$, thus we only need concern ourselves with the
  commutator relations.  Thus we compute:
  \begin{eqnarray*}
    \delta_{x^{p^r}}\cdot \delta_{y^{p^s}} &=& \sum_{a,b\ge 0} \langle m(x)^a m(y)^b, \delta_{x^{p^r}}\tensor \delta_{y^{p^s}}\rangle \delta_{x^a y^b} \\
    &=& \sum_{a,b \ge 0} \langle (x\tensor 1)^a(1\tensor y + x\tensor y)^b,\delta_{x^{p^r}}\tensor \delta_{y^{p^s}}\rangle \delta_{x^a y^b} \\ 
    &=& \sum_{a,b\ge 0} \langle \sum_{b'=0}^b \binom{b}{b'} x^{a+b'}\tensor y^b, \delta_{x^{p^r}}\tensor \delta_{y^{p^s}}\rangle \delta_{x^a y^b} \\ 
    &=& \left\{ 
    \begin{array}{lc}
      \delta_{x^{p^r}y^{p^s}} & r  < s \\
      \delta_{x^{p^r}y^{p^s}} + \delta_{x^{p^r-p^s}y^{p^s}} & r  \ge s.
    \end{array}\right.
  \end{eqnarray*}
  A similar computation shows:
  \begin{displaymath}
    \delta_{y^{p^s}}\cdot \delta_{x^{p^r}} 
    = 
    \sum_{0 \le k \le p^{r-s}}(-1)^k \delta_{x^{p^r - kp^s}y^{p^s}}
  \end{displaymath}
  so that we have as commutator:
  \begin{displaymath}
    \pi_{\canon}(\delta_{x^{p^r}},\delta_{y^{p^s}}) 
    =
    \left\{
    \begin{array}{lc}
      2\delta_{x^{p^r-p^s}y^{p^s}}  - \sum_{2\le k \le p^{r-s}}(-1)^k  \delta_{x^{p^r - kp^s}y^{p^s}} &  r \ge s \\
      0 & r < s.
    \end{array}\right.
  \end{displaymath}
  We note that the calculation is valid even for $p=2$.  The $2$ in the commutator  reflects the fact that the Lie algebra
  of $\mathbb{T}_2$ is abelian for $p=2$.

  Finally, we need to express $\pi_{\canon}$ in multiplicative coordinates.  Taking  $x<y$, we only need
  to apply the following relation recursively:
  \begin{eqnarray*}
    \delta_{x^{p^r}} \delta_{x^my^{p^s}} &=&
    (m_r + 1) \delta_{x^{m+p^r}y^{p^s}} + m_r \delta_{x^m y^{p^s}}\\
    &&
    + \sum_{m\le k \le p^r - p^s} \binom{k}{k-m}\binom{m}{m+p^r -p^s -k} \delta_{x^k y^{p^s}}.
  \end{eqnarray*}
  Here, $m=m_0 + \cdots m_r p^r$ is the $p$-adic expansion of $m$ where we
  assume that $m_r < p-1$.
  \myqed
\end{Example}

\end{document}